\documentclass[12pt,a4paper]{article}
\usepackage[utf8]{inputenc}
\usepackage{amsmath}
\usepackage{amsfonts}
\usepackage{amssymb}
\usepackage{graphicx}
\usepackage{mathabx}
  \usepackage{bbm}
        \usepackage{pdfsync}

\newtheorem{theorem}{Theorem}
\newtheorem{lemma}{Lemma}

\newenvironment{proof}[1][Proof]{\textbf{#1.} }{\ \rule{0.5em}{0.5em}}
\author{Yves Le Jan}
\title{Brownian winding fields } 
\begin{document}
\maketitle
\begin{abstract}

The purpose of the present note is to review and improve the convergence of the renormalized winding fields introduced in \cite{wind1} and \cite{wind2}.
\end{abstract}
%\begin{abstract}
%We investigate the relations between the Poissonnian loop ensembles, their occupation fields, non ramified Galois coverings of a graph, the associated gauge fields, and random Eulerian networks.\\

%Notre étude montre les relations existant entre les ensembles poissoniens de lacets, les champs qu'ils définissent, les circuits euleriens, les revêtements galoisiens des graphes et les champs de jauges associés.
%\end{abstract}

%
%\section{Introduction}
%
In the seminal work of Symanzik \cite{Symanz}, Poisson ensembles of Brownian loops were implicitly used. Since the work of Lawler and Werner \cite{LW} on "loop soups", these ensembles have also been the object of many investigations.\\
Windings of two dimensional random paths have been widely studied. Let us mention the seminal work of Spitzer \cite{spitz} for Brownian paths, and Schramm (\cite{sle}) for SLE. The purpose of the present note is to review and improve the convergence of the renormalized winding fields introduced in \cite{wind1} and \cite{wind2}, using a martingale convergence argument.The result is somewhat reminescent of Gaussian multiplicative chaos ( \cite{Kah}). In the context of Brownian loop ensembles, a different type of renormalization was used to define the occupation field and its powers (see chapter 10 in \cite{stfl}). 
 The method here is related to zeta renormalization used in \cite{lejanbrow}, \cite{lejanbrow2} to study the homology of Brownian loops defined on manifolds.\\

 We consider a bounded open subset of the plane, denoted $D$. We denote by $D_R$ the disc of radius $R$ centered at $0$. For any point $x$ in $D$, let $j_x$ be a uniformizing map mapping $D$ onto $D_1$ and $x$ to $0$ and for $\delta<1$, by $B(x,\delta)$ the pullback of $D_{\delta}$ in $D$.\\
The $\sigma-$finite measure $\mu$ on the set of Brownian loops and the Poisson process of Brownian loops are defined in the same way as Lawler and Werner ''loop soup'' (Cf \cite{LW}).
More precisely, denoting by $dA$ the area measure, $\mu=\int_{x\in X}\int_{0}^{\infty}\frac{1}{t}\mathbb{M}_{t}^{x,x} dt \;dA(x)$ where $\mathbb{M}_{t}^{x,y}$ denotes the distribution of the Brownian bridge in $D$ between $x$ and $y$, multiplied by the heat kernel density $p_t(x,y).$\\
For any positive $\alpha$, the Poisson process of loops of intensity $\alpha \mu$ is denoted $\mathcal{L}_{\alpha}$. If $U$ is an open subset of $D$, we denote by $\mathcal{L}^U_{\alpha}$ the set of loops in $\mathcal{L}_{\alpha}$ contained in $U$.\\
Almost surely, for a given $x$, the loops of $\mathcal{L}_{\alpha}$ do not visit $x$. We denote by $n_x(l)$ the winding number around 0 of the pullback of a loop $l$ in $\mathcal{L}_{\alpha}$. As the Brownian loops, as Brownian paths, have vanishing Lebesgue measure,  $n_x(l)$ is defined almost everywhere in $x$, almost surely.\\
Let $\beta$ denote any $[0, 2\pi)$-valued function defined on $D$. Let $h$ be any bounded function with compact support in $D$. For any $\delta<1$, define
$$W_x ^{\beta_x,\delta,\alpha}=\prod_{l\in \mathcal{L}_{\alpha}\setminus\mathcal{L}^{B(x,\delta)}_{\alpha}}e^{i\beta_x n_x(l)} $$

 The winding field  $W^{\beta, \alpha}(h)$ is defined as follows:

\begin{theorem}

   For $\delta_n$ decreasing to zero, $\int_D h(x)\delta_n^{-\alpha\, a(\beta_x)}W_x ^{\beta_x,\delta_n,\alpha}dA(x)$ is a martingale, with $a(\beta_x)=\frac{\beta_x(2\pi-\beta_x)}{4\pi^2}\leq\frac{1}{4}$. For $\alpha<4$, it converges a.s. and in $L^p$ for all $p\geq 1$ towards a limit denoted $W^{\beta, \alpha}(h)$. 
%The limit, denoted $W^{\beta, \alpha}(h)$.
\end{theorem} 
Remarks: In contrast with Gaussian multiplicative chaos, moment s of all order are defined for any $\alpha<4$.
The question of determining the behaviour of these martingales for $\alpha \geq 4$ seems open. As mentioned in \cite{wind1}, one may also investigate the possibility of finding a characterization of the distribution of the winding field, in terms of conformal field theory. \\

\begin{proof}
 For $0<R \leq \infty$, let $\mathbb{M}_{t}^{R,x,y}$ denote the distribution of the Brownian bridge in $D_R$ multiplied by the heat kernel density,  $\mu^R$ the associated loop measure and $\mathcal{L}_{\alpha}^R$ the corresponding loop ensemble. Up to time change (under which winding indices are invariant), $\mathcal{L}_{\alpha}^1$ is the image of $\mathcal{L}_{\alpha}$ under any uniformizing map.
\begin{lemma}
$ \int_{\mathbb{C}} dA(z)\mathbb{P}_{1}^{\infty,z,z}(n_0=k)=\dfrac{1}{2\pi^2k^2}$
\end{lemma}
This result was established in \cite{Garb}, with reference to \cite{Yor}. Let us outline briefly its proof, for the convenience of the reader:\\
In polar coordinates, a well known consequence of the skew-product decomposition of the Brownian bridge measure is that 
$$\int e^{iu\, n_0(l) }\mathbb{M}_{1}^{\infty,z,z} (dl)=\int e^{iu\int_l d\theta} \mathbb{M}_{1}^{\infty,z,z} (dl)=E(e^{-\frac{u^2}{2}\int_0^1\rho_s^2ds)})q_1(z,z)$$ 
 in which $\rho_s$ denotes a Bessel(0) bridge from $\vert z\vert$ to $\vert z\vert$ and $q_t$ the Bessel(0) transition kernels semigroup.
 It follows from Feynman-Kac formula and Bessel differential equation that this expression equals $e^{-\vert z\vert^2}I_{\vert u\vert}(\vert z\vert)$\\
 As the the Dirac measure at $2\pi n $ is the Fourier transform of $e^{-i 2\pi n u}$ , we get that for any $r>0$
 $$ \mathbb{M}_{1}^{\infty,r,r}(n_0=n)=2e^{-r^2}\int_0^{\infty}I_{\vert u\vert}(r)\cos(2\pi n u) du$$
 From this, as observed by Yor in \cite{Yor}, using the expression of the modified Bessel function $I_{\vert u\vert}$ as a contour integral, we obtain that:
 $$ \mathbb{M}_{1}^{\infty,r,r}(n_0=n)=e^{-r^2}\int_0^{\infty}e^{-r^2\cosh(t)}\left[ \dfrac{2n-1}{r^2+(2n-1)^2\pi^2}-\dfrac{2n+1}{r^2+(2n+1)^2\pi^2}\right]dt$$ 
 
 Hence, integrating with respect to $2\pi r dr$,
 $$ \int_{\mathbb{C}} dA(z)\mathbb{M}_{1}^{\infty,z,z}(n_0=n)=\pi\int_0^{\infty}\frac{dt}{1+\cosh(t)}\left[ \dfrac{2n-1}{r^2+(2n-1)^2\pi^2}-\dfrac{2n+1}{r^2+(2n+1)^2\pi^2}\right]dt.$$
 As observed in \cite{Garb}, the final result follows from a residue calculation yielding telescopic series.
 
 \begin{lemma}
 $\mu^R(l\nsubseteq D_{\delta}, n_0(l)=k) =\frac{1}{2\pi^2k^2}\log(\frac{R}{\delta})$.
 
 \end{lemma}
 To prove this lemma, we use the zeta regularisation method, which, in this context, allows to introduce a $T(l)^s$ factor under $\mu^R$, and let $s$ decrease to zero. ($T(l)$ denoting the loop time length).\\
 $\mu^R(l\nsubseteq D_{\delta},\,  n_0(l)=k)$ is the limit as $s \downarrow 0$ of $\int T(l)^s \mathbbm{1}_{l\nsubseteq D_{\delta} \mathbbm{1}_{n_0(l)=k}}\mu^R(dl) $
 \begin{align*}
 &=\int_0^{\infty} \int_{D_R}\mathbb{M}_{t}^{R,z,z}(n_0=k)dA(z)t^{s-1}dt-\int_0^{\infty} \int _{D_{\delta}}\mathbb{M}_{t}^{\delta,z,z}(n_0=k)dA(z)t^{s-1}dt\\
 &=\int_0^{\infty} \int_{D_R}\mathbb{M}_{t}^{R,z,z}(n_0=k)dA(z) t^{s-1}dt-\int_0^{\infty} \int _{D_{R}}\mathbb{M}_{t(R/ \delta)^2}^{R,z,z}(n_0=k)dA(z)t^{s-1}dt\\
 &=\frac {1-(\delta/R)^{2s}}{s}\int_0^{\infty} \int_{D_R}\mathbb{M}_{t}^{R,z,z}(n_0=k)dA(z) st^{s-1}dt
  \end{align*}
  From lemma 1, for $\eta$ arbitrarily small, we can choose $\epsilon>0$ such that for $u<\epsilon$, $\vert \int_{D_{R/ u}} \mathbb{M}_{1}^{R/u,z,z}(n_0=k)dA(z)-\dfrac{1}{2\pi^2 k^2}\vert<\eta$.\\
  Then $\frac {1-(\delta/R)^{2s}}{s}\int_0^{\epsilon} \int_{D_R}\mathbb{M}_{t}^{R,z,z}(n_0=k)dA(z) st^{s-1}dt=\frac {1-(\delta/R)^{2s}}{s}\int_0^{\epsilon} \int_{D_R/t}\mathbb{M}_{t}^{R/t,z,z}(n_0=k)dA(z) st^{s-1}dt$ is arbitrarily close from $\dfrac{1}{2\pi^2 k^2}\log(\frac{R}{\delta})$ for $\epsilon$ and $s$ small enough.\\
  To prove that $\frac {1-(\delta/R)^{2s}}{s}\int_{\epsilon}^{\infty} \int_{D_R}\mathbb{M}_{t}^{R,z,z}(n_0=k)dA(z) st^{s-1}dt$ converges to zero with $s$, note that $\int_{D_R}\mathbb{M}_{t}^{R}(z,z)(n_0=k)dA(z)\leq \int_{D_R}{P}_{t}^{R}(z,z)dA(z)$, denoting by ${P}_{t}^{R}(x,y)$ the heat kernel on the disc of radius $R$.
  It follows from Weyl asymptotics that this trace can be bounded by $Ce^{-\lambda_0 t}t$, $\lambda_0$ denoting the ground state eigenvalue on $D_R$ and $C$ a positive constant. The result follows as the resulting gamma density converges to zero on $[\epsilon,\infty)$ and this concludes the proof of the second lemma.
  
   \begin{lemma}

 $E(W_x ^{\beta_x,\delta,\alpha})=\delta^{\alpha\, a(\beta_x)}$.
 \end{lemma}
  This result follows by  bounded convergence from lemma 2 and from the Fourier series identity $\sum_1^{\infty}\frac{1}{\pi^2k^2}(1-\cos( k \beta))=\frac{\beta(2\pi-\beta)}{4\pi^2}$  as \\
   \begin{align*}
   E(W_x ^{\beta_x,\delta,\alpha}) &=\lim_{N\rightarrow \infty}E(\prod_{k=-N}^N e^{ik\beta_x\vert \{l\in\mathcal {L}_{\alpha}\,l\nsubseteq B(x,{\delta}), n_x(l)=k\}\vert })\\
  &=\lim_{N\rightarrow \infty}E(\prod_{k=-N}^N e^{ik\beta_x\vert \{l\in\mathcal {L}^1_{\alpha}\,l\nsubseteq D_{\delta}, n_0(l)=k\}\vert })\\
&=\lim_{N\rightarrow \infty}\exp(\alpha \log(\delta)\sum_1^N\frac{1}{\pi^2k^2}(1-\cos( k \beta_x)))
 \end{align*}  
  To complete the proof of the theorem, remark first that it follows from the independence property of a Poisson point process that for $\delta_n$ decreasing to 0, and for any $x$, $ \frac{W_x ^{\beta_x,\delta_n,\alpha}}{E(W_x ^{\beta_x,\delta_n,\alpha})}=\delta_{n}^{-\alpha\, a(\beta_x)}W_x ^{\beta,\delta_n,\alpha}$ is a martingale with independent multiplicative increments. We denote it by $Z_{n,x}^{\beta_x,\alpha}$. Hence, the martingale property of the integral $\int_D h(x) Z_{n,x} dA(x)$  is obvious. To show the convergence, we need a uniform bound on its $L_{2p}$ norm, for any integer $p\geq 1$.\\ 
  Given $2p$ distinct points $x_l$ in a compact $K\subset D$ supporting $h$, for $\delta_{l,n}<\delta_{l,0}=\sup( \{\epsilon, \; B(x_l,\epsilon)\cap B(x_k,\epsilon)=\O\,\, \text{for any } k\neq l\})$ decreasing to zero, all $B(x_l,\delta_{l,0})$ are disjoint and the product $\prod_{l\leq 2p} \delta_{l,n}^{-\alpha\, a(\beta_{x_l})}\,W_{x_l} ^{\beta_{x_l},\delta_{l,n},\alpha}$ is a martingale. 
  Its expectation is bounded by  $\prod_{l\leq 2p} \delta_{l,0}^{-\alpha\, a(\beta_{x_l)}}$.\\
   For some multiplicative constant, $c>0$ depending on the compact support $K$ of $h$, $\delta _{l,0} \leq c\min \{d(x_{l},x_{l'}),\, l'\neq l \}$ for all $l\leq 2p$.
  It follows in particular that $$E(\vert\int_D h(x) Z_{n,x}^{\beta_x,\alpha} dA(x)\vert^2)=E(\int_{D^2} h(x) Z_{n,x}^{\beta_x,\alpha} 
  h(y) Z_{n,y}^{-\beta_y,\alpha} dA(x)dA(y)\leq c^2\Vert h\Vert_{\infty}^2 I $$ with
  $I=\int\int_{D^2}d(x_1,x_2)^{-\alpha/2}dA(x_1)dA(x_2)$, which proves the $L_2$ and a.s. convergence.\\
   More generally, for any integer $p>1$, the $2p$-th moment $E(\vert\int_D h(x) Z_{n,x} dA(x)\vert^{2p})$\\ 
   is bounded by  $ (c\Vert h\Vert_{\infty})^{2p}\int\int_{D^{2p}}\prod_{l\leq 2p}\min_{l'\neq l}d(x_{l'},x_{l})^{- \alpha/4}dA(x_1)...dA(x_{2p})$.
  To see this expression is finite for $\alpha<4$, we will consider only the case $p=2$ as the general proof is similar. The term with highest singularity comes from the case where, up to a permutation, the smallest distances are  $d(x_1,x_2)$ and $d(x_3, x_4)$. Then the integral on that sector of $D^4$ can be bounded by $(I)^2$.
  In the other cases, i.e. when, up to a permutation, the smallest distances are $d(x_1,x_2)$, $d(x_3, x_1)$ and $d(x_4, x_1)$, or $d(x_1,x_2)$, $d(x_3, x_1)$ and $d(x_4, x_2)$, the integral on the corresponding sector can be bounded by $C^2\;I$, with $C=\sup_{x\in K}\int_{D}d(x,y)^{-\alpha/4}dA(y)$. \end{proof}\\
\emph{Remarks:} \\
1) It can be shown that the martingales $Z_{n,x}^{\beta_x,\alpha}$ do not converge, consequently, $W^{\beta, \alpha}(h)$ is a generalized field. The class of test functions $h$ can actually be extended to integrals of delta functions along a smooth curve segment if  $\alpha <4$.\\
2) It follows from theorem 7 in chapter 9 of \cite{stfl}  (see also the Markov property in \cite{wernersemiprob}) that the discrete analogue  of $W^{\beta, \alpha}$ verifies reflection positivity for $\alpha=1$, $2$, or $3$ in case $D$ is invariant under some reflection. This property should extend to the Brownian case.\\

Acknowledgements. I thank Federico Camia and Marci Lis for interesting discussions and the referee for helpful remarks.

\bigskip

\noindent
 NYU-ECNU Institute of Mathematical Sciences at NYU Shanghai. \\   
 D\'epartement de Math\'ematique. Universit\'e Paris-Sud.  Orsay
  
\bigskip
   yves.lejan@math.u-psud.fr \:
   yl57@nyu.edu

\end{document}